\documentclass[12pt]{article}
\usepackage{amssymb,amsmath,epsfig,amscd,eucal}

\newtheorem{thm}{Theorem}[section]
\newtheorem{defn}[thm]{Definition}
\newtheorem{lem}[thm]{Lemma}
\newtheorem{prop}[thm]{Proposition}

\newtheorem{rem}[thm]{Remark}
\newtheorem{cor}[thm]{Corollary}

\newtheorem{letterthm}{Theorem}

\newenvironment{pf}{\par\medskip\noindent{\em Proof. }}{\hfill $\square$\par\medskip}

\newcommand{\Z}{\mathbb{Z}}

\newcommand{\curlyA}{\mathcal{A}}
\newcommand{\curlyB}{\mathcal{B}}

\newcommand{\curlyG}{\mathcal{G}}
\newcommand{\curlyI}{\mathcal{I}}
\newcommand{\curlyL}{\mathcal{L}}
\newcommand{\curlyLR}{\mathcal{LR}}
\newcommand{\curlyS}{\mathcal{S}}

\title{Enumerating limit groups}
\author{Daniel Groves and Henry Wilton}
\date{21st May 2007}

\begin{document}
\maketitle

\begin{abstract}
We prove that the set of limit groups is recursive, answering a question of Delzant.  One ingredient of the proof is the observation that a finitely presented group with local retractions (\emph{\`a la} Long and Reid) is coherent and, furthermore, there exists an algorithm that computes presentations for finitely generated subgroups.  The other main ingredient is the ability to algorithmically calculate centralizers in relatively hyperbolic groups.  Applications include the existence of recognition algorithms for limit groups and free groups.
\end{abstract}

A limit group is a finitely generated, fully residually free group.   Recent research into limit groups has been motivated by their role in the theory of the set of homomorphisms from a finitely presented group to a free group, and in the logic of free groups.  This research has culminated in the independent solutions to Tarski's problems on the elementary theory of free groups by Z.~Sela (see \cite{Se1}, \cite{Se2} \emph{et seq.}) and O.~Kharlampovich and A.~Miasnikov (see \cite{KM98a}, \cite{KM98b} \emph{et seq.}).  Sela's work extends to the elementary theory of hyperbolic groups \cite{Se8}.

We will be entirely concerned with finitely presentable groups.  A class of groups $\curlyG$ is \emph{recursively enumerable} if there exists a Turing machine that outputs a list of presentations for every group $\curlyG$; it is \emph{recursive} if, furthermore, the Turing machine only outputs one presentation from each isomorphism class of $\curlyG$.    T.~Delzant asked if the class of limit groups is recursive \cite[I.13]{SeQ}.

\begin{letterthm}[Corollary \ref{L is recursive}]\label{Theorem A}
The class of limit groups is recursive.
\end{letterthm}

In \cite{BKM} and \cite{DG, DG-Free} it is shown that the isomorphism problem is solvable for the class of limit groups.  Therefore, if the class of limit groups is recursively enumerable it is recursive.  To enumerate limit groups, our approach is to use the structure theory of limit groups developed in \cite{KM98b}.  An equivalent structure theory is described in \cite{Se1}, which could also be used.  Either way, two problems need to be solved.  First, one needs to be able to compute presentations for finitely generated subgroups of limit groups.  We call this property \emph{effective coherence}.  Secondly, one needs to be able to compute centralizers of elements in limit groups.  To solve the second problem we use the relatively hyperbolic structure on limit groups found in \cite{Da} and \cite{Al03}.  Our solution to the first problem relies on local retractions.

D.~Long and A.~Reid \cite{LR2} defined a group to have \emph{local retractions} or \emph{property LR} if every finitely generated subgroup is a retract of a finite-index subgroup.  A finitely presented group with local retractions is coherent.  Furthermore, one can compute presentations for subgroups.

\begin{letterthm}[Theorem \ref{LR is uniformly effectively coherent}]
There exists an algorithm that, given a finite presentation for a group $G$ with local retractions and a finite set of elements $S$, outputs a presentation for the subgroup generated by $S$.
\end{letterthm}

It is a remarkable fact that limit groups are finitely presented.  It was proved in \cite{Wilt4} that limit groups have local retractions.  There is a lengthier proof that limit groups are effectively coherent using the theorem, also proved in \cite{Wilt4}, that iterated centralizer extensions are coset separable with respect to their vertex groups.

As an application of Theorem \ref{Theorem A}, in section \ref{Recognition section} we prove the following theorem.

\begin{letterthm}[Theorem \ref{Recognition theorem}]\label{Theorem C}
There exists an algorithm that, given as input a presentation for a group $G$ and a solution to the word problem in $G$, determines whether or not $G$ is a limit group.
\end{letterthm}

In Corollary \ref{Free group recognition}, we deduce the existence of a similar recognition algorithm for free groups (pointed out to us by Gilbert Levitt).

This paper is the first of a series, in which we intend to prove algorithmic versions of Sela's results.  Specifically, enumerating limit groups will be useful in the algorithmic construction of Makanin--Razborov diagrams over free groups.

\subsection*{Acknowledgements}

The authors would like to thank Zlil Sela for many insightful and generous conversations, and also Fran\c{c}ois Dahmani and Vincent Guirardel for pointing out Corollary \ref{CycPinch} to us.  Thanks also to Gilbert Levitt for drawing Corollary \ref{Free group recognition} to our attention, and to Martin Bridson for explaining the ideas of the paragraph before Theorem \ref{t:findZ}.  The first author was supported in part by NSF Grant DMS-0504251.

\section{Effective coherence}

A finitely generated group is \emph{coherent} if all of its finitely generated subgroups are finitely presented.  We will be interested in the following algorithmic version of coherence.

\begin{defn}
A coherent group $G$ is \emph{effectively coherent} if there exists an algorithm that, given a finite subset $S$ as input, outputs a presentation for the subgroup generated by $S$.

A class $\curlyG$ of coherent groups is \emph{uniformly effectively coherent} if there exists an algorithm that, given as input a presentation of a group $G\in\curlyG$ and a finite set $S$ of elements of $G$, outputs a presentation for the subgroup of $G$ generated by $S$.
\end{defn}

An appealing consequence of this property is that, under mild hypotheses, one can decide if a homomorphism to an effectively coherent group is injective.

\begin{lem}\label{Effective coherence and monomorphisms}
If a group $G$ is effectively coherent then there exists an algorithm that, given a presentation for a group $H$, a solution to the word problem in $H$ and a homomorphism $f:H\to G$, determines whether $f$ is injective.
\end{lem}
\begin{pf}
Given a presentation for the image of $f$ and a solution to the word problem in $H$, it is easy to check whether $f$ has a well-defined inverse and hence is injective.  Therefore, if $G$ is effectively coherent it is easy to check if $f$ is injective.
\end{pf}

\begin{rem} \label{One-sided}
Even without a solution to the word problem in $H$, there exists a Turing machine that will confirm in finite time if the homomorphism $f$ is injective.  Indeed, if $f$ is injective then we know what the inverse to $f$ must be.  By effective coherence, it is possible to compute a presentation for the image of $f$, and the inverse homomorphism exists if and only if the relations for $f(H)$ hold in $H$ (under the supposed inverse map).  Even though the word problem for $H$ may be unsolvable, it is straightforward to enumerate the words which are equal to $1$ in $H$, and if $f$ is a homomorphism then the relations for $f(H)$ (interpreted as words in the generators for $H$) will eventually appear on this list.

However, if the word problem in $H$ is unsolvable then there will in general be no Turing machine which terminates if the map $f$ is not injective, since we will not be able to tell, for example, if the group $H$ is the trivial group.
\end{rem}

Of course, a finitely generated subgroup of an effectively coherent group is effectively coherent.  If $\curlyG$ is a class of groups, denote by $\curlyS(\curlyG)$ the class of finitely generated subgroups of groups in $\curlyG$.  We are interested in effective coherence because it allows the property of being recursively enumerable to pass from $\curlyG$ to $\curlyS(\curlyG)$.  Furthermore, uniform effective coherence also passes to subgroups.

\begin{lem}\label{RE + EC descends to subgroups}
If $\curlyG$ is recursively enumerable and uniformly effectively coherent then $\curlyS(\curlyG)$ is recursively enumerable and uniformly effectively coherent.
\end{lem}
\begin{pf}
Enumerating the presentations of groups $G\in\mathcal{G}$ and finite subsets $S\subset G$, then using uniform effective coherence to compute presentations for $\langle S\rangle$, one enumerates presentations for every group in $\curlyS(\curlyG)$.  So $\curlyS(\curlyG)$ is recursively enumerable.

Given a presentation for a group $G\in\curlyS(\curlyG)$ and a finite subset $S$ of $G$, we can enumerate groups $K\in\curlyG$ and homomorphisms $f: G\to K$ and check whether $f$ is an injection using the Turing machine described in Remark \ref{One-sided}.  Since $G\in\curlyS(\curlyG)$ one will eventually find such an injection $f$.  Using the effective coherence of $K$, one can now compute a presentation for $\langle f(S)\rangle$.  So $\curlyS(\curlyG)$ is uniformly effectively coherent.
\end{pf}

We approach effective coherence through \emph{local retractions}.

\section{Local retractions}

A group $G$ \emph{retracts onto} a subgroup $H$ if the inclusion map $H\hookrightarrow G$ admits a left-inverse $\rho:G\to H$.  The subgroup $H$ is called a \emph{retract} and the map $\rho$ is a \emph{retraction}.  Following \cite{LR2}, a group has \emph{local retractions} if every finitely generated subgroup is a retract of a finite-index subgroup.  This has immediate consequences for coherence.

\begin{lem}\label{Retracts are fp}
If $H$ is a retract of a finitely presented group $G$ then $H$ is finitely presented.
\end{lem}
\begin{pf}
The proof of the lemma is a diagram chase.  Let $\rho:G\to H$ be the retraction.  If $\curlyB$ generates $H$ then, since
$$
G=H\ker\rho
$$
we can add elements from $\ker\rho$ to $\curlyB$ to give a (finite) generating set $\curlyA=\curlyB\cup\curlyA'$ for $G$.  Furthermore, any finite presentation for $G$ can be modified to give a finite presentation with generators of this form.

Denote by $F_X$ the free group on a set $X$.  Let $\rho'$ be the obvious retraction from $F_\curlyA=F_\curlyB \ast F_{\curlyA'}$ to $F_\curlyB$ that kills $F_{\curlyA'}$.  This gives a commutative square
\[ \begin{CD}
{F_\curlyA} @>{p}>> G \\
@V{\rho'}VV    @VV{\rho}V \\
{F_\curlyB} @>{q}>>  {H}\\
\end{CD}\]
where $p$ and $q$ are the natural surjections $F_\curlyA\to G$ and $F_\curlyB\to H$ respectively.  Denote the inclusion $H\hookrightarrow G$ by $i$ and the inclusion $F_\curlyB\hookrightarrow F_\curlyA$ by $i'$.  The lemma follows directly from the claim that $\rho'$ restricts to a retraction $\ker p\to\ker q$.

If $l\in\ker q$ then $p\circ i'(l)=i\circ q(l)=1$ so $i'(l)\in \ker p.$  Likewise, if $k\in\ker p$ then $q\circ\rho'(k)=\rho\circ p(k)=1$ so $\rho'(k)\in\ker q$.  This proves the claim and hence the lemma.
\end{pf}

Since finite-index subgroups of finitely presented groups are finitely presented, coherence for finitely presented groups with local retractions follows immediately.

\begin{prop}
If a finitely presented group $G$ has local retractions then $G$ is coherent.
\end{prop}

Better still, Lemma \ref{Retracts are fp} is effective.

\begin{lem}\label{Computation of presentations of retracts}
Let $G$ be a finitely presented group with solvable word problem.  There is an algorithm that takes as input a finite presentation for $G$ and a collection of words which are the images of the generators under a homomorphism $\rho:G\to G$ that is a retract onto $\rho(G)$, and outputs a presentation for $\rho(G)$.
\end{lem}
\begin{pf}
Applying Tietze transformations, the given generating set for $G$ will eventually be of the form required in the proof of Lemma \ref{Retracts are fp}, namely the union of some generators for $\rho(G)$ and some elements of $\ker\rho$, and since $G$ has solvable word problem we can tell when we have found such a presentation.  By the proof of Lemma \ref{Retracts are fp}, a presentation for $\rho(G)$ is then obtained by eliminating all the generators in $\ker\rho$ from the presentation of $G$.
\end{pf}

By \cite[Theorem 2.4]{LR2}, groups with local retractions are residually finite and hence have (uniformly) solvable word problem.  Let $\curlyLR$ be the class of finitely presented groups with local retractions.

\begin{thm}\label{LR is uniformly effectively coherent}
The class $\curlyLR$ is uniformly effectively coherent.
\end{thm}
\begin{pf}
Given a finite presentation for a group $G\in\curlyLR$ and a finite collection of elements $S\in\curlyG$, we can enumerate all finite-index subgroups $K$ of $G$ using the Reidemeister--Schreier Process (see, for instance, \cite{MKS}).  Since $G \in \curlyLR$, there is a finite-index subgroup $K$ of $G$ so that $\langle S \rangle \subseteq K$ and so that there exists a retraction $\rho : K \to \langle S \rangle$.

We find such a retraction as follows.  In parallel, consider each of the finite-index subgroups of $G$.  Given such a finite-index subgroup $K$, look for the elements of $S$ as words in the generators for $K$.  Suppose we have found a finite-index subgroup $K$ so that $\langle S \rangle \subseteq K$, and a finite presentation $\langle X \mid R(X) \rangle$ of $K$, with $S = \{ s_1(X), \ldots , s_n(X) \}$ written as words in $X^\pm$.  Now search for a collection of words $Y$ in $S^\pm$ with a bijection $\rho : X \to Y$ so that each of the relations of the form $R(Y)$ holds and so that for each $i$ we have $s_i(Y) = s_i(X)$.  Then the map $\rho$ extends to a retraction $\rho : K \to \langle S\rangle$.  Since there is a retraction, we will eventually find such a $K$ and $Y$.

The algorithm of Lemma \ref{Computation of presentations of retracts} now computes a presentation for $\langle S\rangle$.
\end{pf}

\section{Enumerating $\curlyI$ and $\curlyL$}

The class of \emph{iterated extensions of centralizers} is defined inductively.  If $G$ is a group, $g\in G$ and $Z(g)$ is the centralizer of $g$ then an amalgamated free product
$$
G'=G*_{Z(g)} (Z(g)\times \Z^n)
$$
is said to be obtained from $G$ by \emph{extension of centralizers}.

\begin{defn}
The class $\curlyI$ of \emph{iterated extensions of centralizers} is the smallest class of groups containing all finitely generated free groups and closed under extension of centralizers.  The class of \emph{limit groups} is defined to be
$$
\curlyL=\curlyS(\curlyI),
$$
the class of finitely generated subgroups of iterated extensions of centralizers.
\end{defn}
The usual definition of limit groups is as finitely generated fully residually free groups.

\begin{defn}
A group $G$ is \emph{fully residually free} if, for every finite subset $X\subset G\smallsetminus 1$, there exists a homomorphism to a free group $G\to F$ such that $1\notin f(X)$.
\end{defn}
A finitely generated group is fully residually free if and only if it is in $\curlyL$, by a theorem of \cite{KM98b}.  Fully residually free groups are residually finite (since free groups are) and so have solvable word problem.  Using the fact that limit groups are fully residually free, the following fact is well known and easy to prove.

\begin{lem}
If $G$ is a limit group and $g\in G$ then $Z(g)$ is a free abelian group.
\end{lem}

By Theorem B of \cite{Wilt4}, limit groups have local retractions.  It is clear that all groups in $\curlyI$ are finitely presented.

\begin{cor}\label{I is uniformly coherent}
The class $\curlyI$ is uniformly effectively coherent.
\end{cor}

By Lemma \ref{RE + EC descends to subgroups}, to enumerate limit groups it remains only to enumerate $\curlyI$.  The crucial step is the ability to calculate centralizers.

For this we use the relatively hyperbolic structure of limit groups (found independently by E.~Alibegovi\'{c} \cite{Al03} and F.~Dahmani \cite{Da}).  See \cite{Farb} for an introduction to relatively hyperbolic groups (where in Farb's language we mean `relatively hyperbolic with BCP').  Limit groups are torsion-free and hyperbolic relative to a finite collection of maximal noncyclic abelian subgroups.  Dahmani \cite{Da2} provides an algorithm which takes as input a finite presentation of such a relatively hyperbolic group and outputs a basis for a representative of each conjugacy class of noncyclic maximal abelian subgroup (Dahmani's algorithm takes as input an arbitrary finite presentation, and does not need to be given the `relatively hyperbolic structure' of the group).

Another important tool will be the universal theory of a group.   The \emph{elementary theory} of a group $G$ is the set of all sentences in first-order predicate logic (possibly with constants) that hold in $G$.  For example, $G$ is abelian if and only if the sentence
\[
\forall x,y\in G\  [x,y]=1
\]
is in the elementary theory of $G$.  A \emph{universal sentence} is a sentence in the elementary theory with a single universal quantifier.  The \emph{universal theory} of $G$ is the set of universal sentences in the elementary theory of $G$.  Deciding the truth of universal sentences is equivalent to deciding whether finite systems of equations and inequations (with constants) have solutions.

In \cite{Makanin}, Makanin proved that the universal theory of a free group $F$ is \emph{decidable}---that is, there exists an algorithm that, given as input a universal sentence, determines whether or not it lies in the universal theory of $F$.   The universal theory of torsion-free relatively hyperbolic groups with abelian parabolic subgroups is also decidable, by another algorithm of Dahmani \cite{Da3} (again the input is any finite presentation for the group, along with the universal sentence).

There is an alternative approach to calculating centralizers using biautomatic structures.  It follows from work of Rebbechi \cite{Reb} that limit groups are biautomatic, and the algorithm for finding automatic structures described in \cite{ECHLPT} can be adapted to find biautomatic structures \cite{BR}.  In particular, one can calculate the fellow-traveller constant of the bicombing.  Using the ideas of \cite{Br01}, it is then easy to compute a presentation for the centralizer of an arbitrary finite subset.

\begin{thm} \label{t:findZ}
There exists an algorithm that, given as input a presentation for a group $G\in \curlyI$ and an element $g\in G$, outputs a minimal set of generators for $Z(g)$.
\end{thm}
\begin{pf}
Apply Dahmani's algorithm from \cite{Da2} to find a basis for a representative of each conjugacy class of maximal noncyclic abelian subgroup.

Let $g \in G$.  There are two cases to consider: either $g$ is {\em parabolic} (which means conjugate into a noncyclic abelian subgroup) or else $g$ is {\em hyperbolic} (which means $g$ is not parabolic).

It is possible to decide whether or not $g$ is parabolic.  This is because the universal theory of $G$ is decidable \cite{Da3}.  The element is parabolic if and only if there exists an element $h \in G$ so that $hgh^{-1}$ commutes with each element of one of the above bases for the noncyclic abelian subgroups.  This is a finite system of equations  over $G$, which we can determine the truth of by Dahmani's algorithm from \cite{Da3}.

If $g$ is parabolic, then we will find such an element $h$, and the conjugates by $h^{-1}$ of the basis for the maximal noncyclic abelian subgroup generates the centralizer of $g$.  In this case we have found a minimal generating set for $Z(g)$.

If $g$ is hyperbolic then its centralizer is generated by a maximal root of $g$.  According to D.~Osin \cite[Theorem 1.16.(3)]{Osin}, it is possible to algorithmically extract roots from hyperbolic elements of $G$.  On the face of it, Osin's algorithm needs to be given as input the relatively hyperbolic structure of the group.  However, Dahmani's algorithm from \cite{Da2} will find this structure, so we can make Osin's algorithm take only the finite presentation as input.  Therefore, if $g$ is hyperbolic we can find a maximal root of $g$, and this maximal root is a minimal generating set for $Z(g)$.
\end{pf}

\begin{cor}
The set $\curlyI$ is recursively enumerable.
\end{cor}

Combining this with Theorem \ref{I is uniformly coherent} it follows that the set of limit groups $\curlyL$ is recursively enumerable, by Lemma \ref{RE + EC descends to subgroups}.

\begin{cor}
The set of limit groups $\curlyL$ is recursively enumerable and uniformly effectively coherent.
\end{cor}

The results of \cite{BKM} (see also \cite{DG, DG-Free}) show that limit groups have solvable isomorphism problem.  Hence we can improve recursively enumerable to recursive: we can ensure that the list produced includes at most one presentation for each isomorphism class of limit groups.

\begin{cor}\label{L is recursive}
The set of limit groups $\curlyL$ is recursive.
\end{cor}

On the other hand, by systematically applying Tietze transformations, it is possible to effectively list all of the finite presentations of all limit groups.  We give a simple application to recognition algorithms here.

\section{Recognition algorithms}\label{Recognition section}

\begin{thm}\label{Recognition theorem}
There exists an algorithm that, given as input a presentation for a group $G$ and a solution to the word problem in $G$, determines whether or not $G$ is a limit group.
\end{thm}
\begin{pf}
Let $\mathcal{P} = \langle X \mid R \rangle$ be the finite presentation defining $G$.  We have already noted that it is possible to enumerate all finite presentations of limit groups.  Thus if $G$ is a limit group then $\mathcal P$ will eventually appear on this list.

Suppose then that $G$ is not a limit group.  Then $G$ is not fully residually free, so there is a finite set $\{ g_1, \ldots , g_r \}$ of nontrivial elements of $G$ so that for any homomorphism $\phi$ from $G$ to a free group  $F$, at least one of the $g_i$ is in $\ker(\phi)$.  This property of $G$ can easily be translated into a system of equations and inequations over $F$ as follows.  Consider both the elements of $R$ and each $g_i$ as a word in $X^\pm$, and write $R = \{ r_1 , \ldots , r_k \}$.  Then the following sentence encodes the fact that at least one of $\{ g_1, \ldots , g_r \}$ is in the kernel of any homomorphism from $G$ to $F$:
\begin{equation} \label{Sentence}
 \forall X \subset F\  \big(r_1(X) = 1 \wedge \cdots \wedge r_k(X) = 1 \big) \Rightarrow
\big(g_1(X) = 1 \vee \cdots \vee g_r(X) = 1\big).
\end{equation}
By Makanin's algorithm \cite{Makanin}, it is possible to decide whether or not universal sentences are true in a free group.  Enumerate finite sets of nontrivial elements of $G$ (the solution to the word problem allows us to know that the elements are nontrivial).  Now, for each such finite set $\{ r_1, \ldots , r_k \}$, decide whether the sentence \eqref{Sentence} is true or not.  If $G$ is not a limit group, we will eventually find a finite set for which \eqref{Sentence} is true.
\end{pf}

Of course, one cannot recognize limit groups amongst arbitrary finitely presented groups.

A {\em cyclically pinched group} is an amalgamated free product of two free groups with cyclic amalgamated subgroup.  Some, but not all, of these groups are limit groups.  In \cite[I.3]{SeQ}, Sela asks for necessary or sufficient conditions for a cyclically pinched group to be a limit group.  We do not have an answer to this question.  However, Theorem \ref{Recognition theorem} implies that at least the question has an answer.
The following result was pointed to us by Fran\c{c}ois Dahmani and Vincent Guirardel
(its proof contains the core of the proof of Theorem \ref{Recognition theorem}).

\begin{cor} \label{CycPinch}
There is an algorithm that takes as input a finite presentation of a cyclically pinched group and decides whether or not the defined group is a limit group.
\end{cor}

It does not matter whether the input presentation exhibits the cyclically pinched nature of the group, since by applying some finite number of Tietze transformations it is possible to find such a presentation.  Once such a presentation is found, there is an explicit solution to the word problem. Therefore Corollary \ref{CycPinch} follows immediately from Theorem \ref{Recognition theorem}.

As remarked above, limit groups are torsion-free and hyperbolic relative to their maximal abelian subgroups.  There is an algorithm to distinguish free groups among such relatively hyperbolic groups; indeed, it is proved in \cite[Theorem 1.4]{DG-Free} that there exists an algorithm that computes the Grushko decomposition from a presentation of such a group.  Combining this with Theorem \ref{Recognition theorem}, we obtain a similar recognition algorithm for free groups.  This corollary was pointed out to us by Gilbert Levitt.

\begin{cor}\label{Free group recognition}
There exists an algorithm that, given as input a presentation for a group $G$ and a solution to the word problem in $G$, determines whether or not $G$ is free.
\end{cor}

One can also deduce a similar result for surfaces.  In \cite[Theorem D]{DG} it is shown that there exists an algorithm that computes a JSJ decomposition for a torsion-free, freely indecomposable group that is hyperbolic relative to its maximal abelian subgroups.  In particular, combining this with the algorithm from \cite{DG-Free}, one can decide whether or not a limit group is a surface group.  It follows as before that there exists an algorithm that, given as input a presentation for a group $G$ and a solution to the word problem in $G$, determines whether or not $G$ is a (fully) residually free surface group.  (The only surface groups that are not residually free are the fundamental groups of the non-orientable surfaces of Euler characteristic 1, 0 and -1.)

\bibliographystyle{plain}

\begin{thebibliography}{10}

\bibitem{Al03}
Emina Alibegovi{\'c}.
\newblock A combination theorem for relatively hyperbolic groups.
\newblock {\em Bull. London Math. Soc.}, 37(3):459--466, 2005.

\bibitem{Br01}
Martin~R. Bridson.
\newblock On the subgroups of semihyperbolic groups.
\newblock In {\em Essays on geometry and related topics, Vol. 1, 2}, volume~38
  of {\em Monogr. Enseign. Math.}, pages 85--111. Enseignement Math., Geneva,
  2001.

\bibitem{BR}
Martin~R. Bridson and Lawrence~D. Reeves.
\newblock On the algorithmic construction of classifying spaces and the
  isomorphism problem for biautomatic groups.
\newblock Preprint, 2007.

\bibitem{BKM}
I.~Bumagin, O.~Kharlampovich, and A.~Miasnikov.
\newblock Isomorphism problem for finitely generated fully residually free
  groups.
\newblock {\em J. Pure and Applied Algebra}, 208(3):961--977, 2007.

\bibitem{Da3}
Fran\c{c}ois Dahmani.
\newblock Existential questions in (relatively) hyperbolic groups.
\newblock Preprint, 2006.

\bibitem{Da2}
Fran\c{c}ois Dahmani.
\newblock Finding relatively hyperbolic structures.
\newblock Preprint, 2006.

\bibitem{DG-Free}
Fran\c{c}ois Dahmani and Daniel Groves.
\newblock Detecting free splittings in relatively hyperbolic groups.
\newblock TAMS, to appear.

\bibitem{DG}
Fran\c{c}ois Dahmani and Daniel Groves.
\newblock The isomorphism problem for toral relatively hyperbolic groups.
\newblock Preprint, 2005.

\bibitem{Da}
Fran{\c{c}}ois Dahmani.
\newblock Combination of convergence groups.
\newblock {\em Geom. Topol.}, 7:933--963 (electronic), 2003.

\bibitem{ECHLPT}
David B.~A. Epstein, James~W. Cannon, Derek~F. Holt, Silvio V.~F. Levy,
  Michael~S. Paterson, and William~P. Thurston.
\newblock {\em Word processing in groups}.
\newblock Jones and Bartlett Publishers, Boston, MA, 1992.

\bibitem{Farb}
B.~Farb.
\newblock Relatively hyperbolic groups.
\newblock {\em Geom. Funct. Anal.}, 8(5):810--840, 1998.

\bibitem{KM98a}
O.~Kharlampovich and A.~Miasnikov.
\newblock Irreducible affine varieties over a free group. {I}. {I}rreducibility
  of quadratic equations and {N}ullstellensatz.
\newblock {\em J. Algebra}, 200(2):472--516, 1998.

\bibitem{KM98b}
O.~Kharlampovich and A.~Miasnikov.
\newblock Irreducible affine varieties over a free group. {II}. {S}ystems in
  triangular quasi-quadratic form and description of residually free groups.
\newblock {\em J. Algebra}, 200(2):517--570, 1998.

\bibitem{LR2}
D.~D. Long and A.~W. Reid.
\newblock Subgroup separability and virtual retractions of groups.
\newblock \emph{Topology}, to appear, 2006.

\bibitem{MKS}
Wilhelm Magnus, Abraham Karrass, and Donald Solitar.
\newblock {\em Combinatorial group theory: {P}resentations of groups in terms
  of generators and relations}.
\newblock Interscience Publishers [John Wiley \& Sons, Inc.], New
  York-London-Sydney, 1966.

\bibitem{Makanin}
G.~S. Makanin.
\newblock Decidability of the universal and positive theories of a free group.
\newblock {\em Izv. Akad. Nauk SSSR Ser. Mat.}, 48(4):735--749, 1984.

\bibitem{Osin}
Denis~V. Osin.
\newblock Relatively hyperbolic groups: intrinsic geometry, algebraic
  properties, and algorithmic problems.
\newblock {\em Mem. Amer. Math. Soc.}, 179(843):vi+100, 2006.

\bibitem{Reb}
D.~Y. Rebbechi.
\newblock Algorithmic properties of relatively hyperbolic groups.
\newblock Thesis, 2003. ArXiv:~math/0302245.

\bibitem{Se8}
Z.~Sela.
\newblock Diophantine geometry over groups {VIII}: The elementary theory of a
  hyperbolic group.
\newblock Preprint.

\bibitem{SeQ}
Zlil Sela.
\newblock Diophantine geometry over groups: a list of research problems.
\newblock \texttt{http://www.ma.huji.ac.il/\textasciitilde zlil/problems.dvi}.

\bibitem{Se1}
Zlil Sela.
\newblock Diophantine geometry over groups. {I}. {M}akanin-{R}azborov diagrams.
\newblock {\em Publ. Math. Inst. Hautes \'Etudes Sci.}, 93:31--105, 2001.

\bibitem{Se2}
Zlil Sela.
\newblock Diophantine geometry over groups. {II}. {C}ompletions, closures and
  formal solutions.
\newblock {\em Israel J. Math.}, 134:173--254, 2003.

\bibitem{Wilt4}
Henry Wilton.
\newblock Hall's {T}heorem for limit groups.
\newblock Preprint, 2006.

\end{thebibliography}

\noindent
\textsc{Daniel Groves, Mathematics 253-37, California Institute of Technology, Pasadena, CA 91125}\\
\emph{E-mail:} \texttt{groves@caltech.edu}\\\\
\textsc{Henry Wilton, Department of Mathematics, 1 University Station C1200, Austin, TX 78712-0257}\\
\emph{E-mail:} \texttt{henry.wilton@math.utexas.edu}\\

\end{document}